\documentclass[12pt]{amsart}
\usepackage{amsmath}
\usepackage{amssymb}
\usepackage{latexsym}

\newcommand{\Zint}{\mathbb {Z}}    
     
\newcommand{\Rea}{\mathbb {R}}      
\newcommand{\Cplx}{\mathbb {C}}     

\newcommand{\halmos}{\rule{5pt}{5pt}}

\numberwithin{equation}{section}

\newtheorem{prop}{\bf Proposition}[section]
\newtheorem{thm}[prop]{\bf Theorem}

\newenvironment{rmk}{\noindent{\bf Remark}\hskip 5pt}{\hfill{$\Box$}}

\begin{document}

\title[Inozemtsev models]
{Quasi-exact solvability of Inozemtsev models}
\author{Kouichi Takemura}
\address{Department of Mathematical Sciences, Yokohama City University, 22-2 Seto, Kanazawa-ku, Yokohama 236-0027, Japan.}
\email{takemura@yokohama-cu.ac.jp}

\begin{abstract}
Finite-dimensional spaces which are invariant under the action of the Hamiltonian of the $BC_N$ Inozemtsev model are introduced, and it is shown that commuting operators of conserved quantities also preserve the finite-dimensional spaces.

The relationship between the finite-dimensional spaces of the $BC_N$ Inozemtsev models and the theta-type invariant spaces of the $BC_N$ Ruijsenaars models is clarified. The degeneration of the $BC_N$ Inozemtsev models and the correspondence of their invariant spaces are considered.
\end{abstract}

\maketitle

\section{Introduction}

In \cite{Ino}, Inozemtsev proposed a $N$-particle quantum mechanics model, which is called the $BC_N$ Inozemtsev model.
This is a generalization of the Calogero-Moser-Sutherland model or the Olshanetsky-Perelomov model \cite{OP}.

The $BC_N$ Inozemtsev model is quantum completely integrable.
Here, the quantum complete integrability means that there exists $N$ algebraically independent mutually commuting operators (higher commuting Hamiltonians) which commute with the Hamiltonian of the model. It is a quantum version of Liouville's integrability.
For the $BC_N$ Inozemtsev model, Oshima \cite{O} described the commuting operators explicitly.
Note that the $BC_N$ Inozemtsev model is a universal completely integrable model of quantum mechanics with $B_N$ symmetry, which follows from the classification due to Ochiai, Oshima and Sekiguchi \cite{OOS,OS}.

On the other hand, Finkel, Gomez-Ullate, Gonzalez-Lopez, Rodriguez and Zhdanov studied quasi-exactly solvable models in \cite{FGGRZ1,FGGRZ2}. They found several quasi-exactly solvable many-body systems.
Although they did not use the phrase ``$BC_N$ Inozemtsev model'', they essentially found that the $BC_N$ Inozemtsev model is quasi-exactly solvable, i.e. the Hamiltonian of the $BC_N$ Inozemtsev model preserves some finite-dimensional space which is spanned by some symmetric ``monomials''. 

In this paper, we link the quasi-exact solvability with the quantum complete integrability. More precisely, we show that the commuting operators (higher Hamiltonians) of the $BC_N$ Inozemtsev model also preserve the finite-dimensional space, which has been appeared in the context of quasi-exact solvability.

On the finite-dimensional space, joint eigenvalues and eigenfunctions of the commuting operators are determined by algebraic calculations. In this sense, the model would be solved partially. Note that the phrase ``quasi-exact solvability'' is used in these situations (see \cite{Tur}).

The spectral problem of quantum mechanics is generally considered in a Hilbert space, and the Hilbert space is often taken as a square-integrable space ($L^2$ space). Therefore it would be important to consider the relationship between the Hilbert space ($L^2$ space) of the $BC_N$ Inozemtsev model and the finite-dimensional space which appears in the context of the quasi-exact solvability. In this paper, we determine the condition for that the finite-dimensional space lies in the $L^2$ space.

There are other models which are concerned with the quasi-exactly solvability.
Ruijsenaars-type models are introduced for arbitrary root systems including the $BC_N$ cases, which are difference analogues of Inozemtsev (or Calogero-Moser-Sutherland) models (see \cite{Rui,vD,HK}). Hikami and Komori \cite{HK,Kom1,Kom2,KH} finally constructed higher commuting operators using algebra, and found an invariant subspace spanned by theta functions.
In \cite{ST}, Sasaki and Takasaki considered degenerate Inozemtsev models and their quasi-exact solvability.

In this paper, correspondence between the $BC_N$ Ruijsenaars model and the $BC_N$ Inozemtsev model is considered. We will observe that the invariant subspace spanned by theta functions for the $BC_N$ Ruijsenaars model corresponds to the invariant space related to the quasi-exact solvability for the $BC_N$ Inozemtsev model. 

To obtain the degenerate $BC_N$ Inozemtsev model from the (elliptic) $BC_N$ Inozemtsev model, a certain trigonometric limit is considered. It is shown that the finite-dimensional invariant spaces for the $BC_N$ Inozemtsev model tend to the invariant spaces of degenerate $BC_N$ Inozemtsev model which were introduced by Sasaki and Takasaki. It would be important to consider degeneration of models, because it would be helpful in understanding several integrable models and the relationship among them.

This paper is organized as follows. 
In section \ref{sec:FGGRZ}, the $BC_N$ Inozemtsev model and its finite-dimensional invariant spaces are introduced.
In section \ref{sec:comm}, higher commuting operators of the $BC_N$ Inozemtsev model are introduced. In section \ref{sec:commQES}, it is shown that higher commuting operators also preserve the finite-dimensional invariant spaces. 
In section \ref{sec:L2}, the relationship between the finite-dimensional invariant spaces and the $L^2$ space is considered. 
If the coupling constants are integers, the model may have some special features. In section \ref{sec:int}, we consider this case.
In section \ref{sec:Ruij}, correspondence between the theta-type invariant spaces for the $BC_N$ Ruijsenaars model and the invariant spaces which is related to the quasi-exact solvability for the $BC_N$ Inozemtsev model are investigated.
In section \ref{sec:deg}, we consider the degeneration of the $BC_N$ Inozemtsev models and the limit of invariant spaces.

{\bf Acknowledgment}
The author would like to thank Professors ~Yasushi Komori, ~Toshio Oshima,  ~Ryu Sasaki, and ~Kanehisa Takasaki for fruitful discussions. He is partially supported by the Grant-in-Aid for Scientific Research (No. 13740021) from the Japan Society for the Promotion of Science.

\section{$BC_N$ Inozemtsev model and its invariant space} \label{sec:FGGRZ}

The $BC_N$ Inozemtsev model is a system of quantum mechanics with $N$-particles whose Hamiltonian is given by
\begin{align}
& H=-\sum_{j=1}^N\frac{\partial ^2}{\partial x_j^2}+2l(l+1)\sum_{1\leq j<k\leq N} \left( \wp (x_j-x_k) +\wp (x_j +x_k) \right) \label{InoHam} \\
& \; \; \; \; \; + \sum_{j=1}^N \sum _{i=0}^3 l_i(l_i+1) \wp(x_j +\omega_i), \nonumber
\end{align}
where $\wp (x)$ is the Weierstrass' $\wp$-function with periods $(1,\tau )$ (see (\ref{def:wp})), $\omega _0=0$, $\omega_1=\frac{1}{2}$, $\omega_2=-\frac{\tau+1}{2}$, $\omega_3=\frac{\tau}{2}$ are half periods, and $l$ and $l_i$ $(i=0,1,2,3)$ are coupling constants.

Let $a$, $b_i$ $(i=0,1,2,3)$ be numbers which satisfy $a\in \{ -l, l+1 \}$ and $b_i \in \{-\frac{l_i}{2}, \frac{l_i+1}{2}\}$ $(i=0,1,2,3)$.
Set $z_j=\wp (x_j)$ $(1\leq j\leq N)$ and 
\begin{align}
& \Phi(z)=\prod_{1\leq j<k\leq N} (z_j-z_k)^a \prod_{j=1}^N \prod _{i=1}^3 (z_j-e_i)^{b_i}, \label{def:Phi}\\
& \widehat{H}= \Phi(z)^{-1} \circ H \circ  \Phi(z), \nonumber
\end{align}
where $e_i=\wp(\omega_i)$, $(i=1,2,3)$.

By applying formulae (\ref{wpformula}), it is shown directly that the operator $\widehat{H}$ admits the following expression:
\begin{align}
&\widehat{H}= -\left( \sum_{j=1}^N 4(z_j-e_1)(z_j-e_2)(z_j-e_3) \left( \frac{\partial ^2}{\partial z_j^2} \right. \right. \label{gtHam} \\
&+ \left. \left. \left( \sum_{k\neq j}\frac{2a}{z_j-z_k}+\frac{2b_1+\frac{1}{2}}{z_j-e_1}+\frac{2b_2+\frac{1}{2}}{z_j-e_2}+\frac{2b_3+\frac{1}{2}}{z_j-e_3}\right)\frac{\partial}{\partial z_j} \right) \right) \nonumber \\
&-4\left((N-1)a-b_0+b_1+b_2+b_3+\frac{1}{2}\right)\left((N-1)a+b_0+b_1+b_2+b_3\right)\left( \sum_{j=1}^Nz_j \right) \nonumber \\
&+4N((b_1+b_2)^2e_3+(b_1+b_3)^2e_2+(b_2+b_3)^2e_1) \nonumber \\
& -4N(N-1)a(e_1b_1+e_2b_2+e_3b_3).\nonumber
\end{align}

\begin{prop} \label{prop:Hinv} 
(i) Let $P^{\mbox{\rm \scriptsize sym}}$ be the space of symmetric polynomials in variables $z_1, z_2, \dots ,z_N$. Then $\widehat{H} \cdot P^{\mbox{\rm \scriptsize sym}} \subset P^{\mbox{\rm \scriptsize sym}}$.\\
(ii) Let $a$, $b_i$ $(i=0,1,2,3)$ be numbers which satisfy $a\in \{ -l, l+1 \}$ and $b_i \in \{-\frac{l_i}{2}, \frac{l_i+1}{2}\}$ $(i=0,1,2,3)$. Assume that $d=-((N-1)a+b_0+b_1+b_2+b_3)$ is a non-negative integer. Let $V_d$ be the vector space spanned by monomials $z_1^{m_1}z_2^{m_2}\dots z_N^{m_N} $ such that $m_i \in \{ 0,1, \dots ,d\}$ for all $i$, and $V_d^{\mbox{\scriptsize sym}}=V_d \cap P^{\mbox{\rm \scriptsize sym}}$. Then $\widehat{H}  \cdot V_d^{\mbox{\rm \scriptsize sym}} \subset V_d^{\mbox{\rm \scriptsize sym}}$.
\end{prop}
\begin{proof}
Let $f \in P^{\mbox{\scriptsize sym}}$. From (\ref{gtHam}), the function $\widehat{H} f$ is a symmetric rational function and it may have poles only along $z_j-z_k=0$ of degree at most $1$. If $\widehat{H} f$ has a pole along $z_j-z_k=0$ $(j\neq k)$ of degree $1$, it contradicts to the symmetry of $\widehat{H} f$ on the variables $z_j$ and $z_k$. Hence the function $\widehat{H} f$ does not have poles and we obtain (i).

Let $\lambda _1, \dots , \lambda _N$ be non-negative integers such that $\lambda _1\geq \lambda _2 \geq \dots \geq \lambda _N \geq 0$.

Set $P(z)=4(z-e_1)(z-e_2)(z-e_3)$ and $L=\lambda _1$.
Let $l$ and $l'$ be integers such that $0\leq l,l' \leq L$. Then
\begin{align}
& \left( \frac{P(z_j)}{(z_j-z_k)}\frac{\partial }{\partial z_j}+\frac{P(z_k)}{(z_k-z_j)}\frac{\partial }{\partial z_k} \right) (z_j^{l}z_k^{l'}+z_j^{l'}z_k^{l}) \label{eq:part} \\
& = 4L \delta _{l,L} (z_j^{L +1}z_k^{l'}+z_j^{l'}z_k^{L +1}) +4L \delta _{l',L } (z_j^{l}z_k^{L +1}+z_j^{L +1}z_k^{l}) +(\#_1) ,\nonumber
\end{align}
where $\delta _{l,L}$ is the Kronecker's delta and the term $(\# _1)$ is a linear combination of monomials $z_j^t z_k^{t'}$ such that $0\leq t,t'\leq L$.

Let $\lambda = (\lambda _1, \lambda _2, \dots ,\lambda _N)$ and $\lambda ^+= (\lambda _1+1, \lambda _2, \dots ,\lambda _N)$. By summing up the equality (\ref{eq:part}) on $j$ and $k$ $(1\leq j< k\leq N)$, we obtain
\begin{align}
& \sum_{j=1}^{N} \sum_{k \neq j} \frac{P(z_j)}{(z_j-z_k)}\frac{\partial }{\partial z_j} m_{\lambda } = 4 (N-1)L m_{\lambda ^+}+ (\#_2),
\end{align}
where $m_{\mu }= \sum_{(m_1, \dots ,m_N) \in {\frak S}_N \cdot \mu } z_1^{m_1}z_2^{m_2}\dots z_N^{m_N}$ for $\mu =(\mu _1, \dots ,\mu _N)$, ${\frak S}_N$ is the symmetric group of $N$ letters and the term $(\# _2)$ is a linear combination of symmetric monomials $m_{(\mu _1, \dots ,\mu _N)}$ such that $0\leq \mu _N \leq  \dots \leq \mu_1 \leq L$.

Hence we obtain
\begin{align}
& \widehat{H} m_{\lambda } = -4(L-d-2b_0+\frac{1}{2})(L-d) m_{\lambda ^+} + (\#_3),
\label{Hactml}
\end{align}
where the term $(\# _3)$ is a linear combination of symmetric monomials $m_{(\mu _1, \dots ,\mu _N)}$ such that $0\leq \mu _N \leq  \dots \leq \mu_1 \leq L$.

If $L\leq d$ then all elements in the term $(\# _3)$ lie in $V_d^{\mbox{\scriptsize sym}}$.
If $L< d$ then $m_{\lambda ^+} \in V_d^{\mbox{\scriptsize sym}}$, and if $L=d$ then the coefficient of $m_{\lambda ^+}$ on the right-hand side of (\ref{Hactml}) vanishes. Hence if $L\leq d$ then $\widehat{H} m_{\lambda } \in V_d^{\mbox{\scriptsize sym}}$.

Therefore we obtain (ii).
\end{proof}
\begin{rmk}
Proposition \ref{prop:Hinv} was essentially obtained in \cite{FGGRZ2} by a different method.
\end{rmk}

\section{Commuting operators and invariant subspaces}

\subsection{Commuting operators} \label{sec:comm} $ $

It is known that the $BC_N$ Inozemtsev model is completely integrable, i.e. there exists $N$ algebraically independent mutually commuting operators (higher commuting Hamiltonians) which commute with the Hamiltonian.
In \cite{O}, Oshima gave explicit forms of the commuting operators.
Now we pick up some results obtained by Oshima.

Let $W(B_N)$ be the Weyl group of type $B_N$, i.e. the group of the coordinate transformations 
\begin{equation}
(x_1, \dots , x_N) \mapsto (\epsilon _1 x_{\sigma (1)}, \dots , \epsilon _N x_{\sigma (N)}) \label{BNtrans}
\end{equation}
of $\Rea ^N$, where $\sigma \in {\mathfrak S}_N$ (the symmetric group) and $\epsilon _1=\pm 1, \dots ,\epsilon _N=\pm 1 $. Let $W(D_N)$ be a subgroup of $W(B_N)$ which consist of transformations (\ref{BNtrans}) with a condition $\prod _{j=1}^N \epsilon _j =1$. For $w \in W(B_N)$ we define $\epsilon (w)=\left\{ \begin{array}{ll} 1 & w \in W(D_N) \\ -1 & w \not\in W(D_N)\end{array} \right. $.

Let us consider the operators which commute with the Hamiltonian of the $BC_N$ Inozemtsev model (\ref{InoHam}).
For $i=0,1,2,3$, we set
\begin{align}
& S_{\{1,\dots ,k\}}= \sum_{w\in W(B_k)}w\left( \wp(x_1-x_2)\wp (x_2-x_3) \dots \wp (x_{k-1}-x_{k})\right), \\
& S_{\{1,\dots ,k\}} ^{\langle i \rangle } = \sum_{w\in W(B_k)}w\left( \wp (x_1+\omega _i) \wp(x_1-x_2)\wp (x_2-x_3) \dots \wp (x_{k-1}-x_{k})\right), \nonumber \\
& T^o_{\{1,\dots ,k\}} = \sum_{I_1 \amalg \dots \amalg I_{\mu} = \{1,\dots ,k \}} (-1)^{\mu -1} (\mu -1)! S_{I_1} \dots S_{I_{\mu }},  \nonumber\\
& T^{o,\langle i \rangle } _{\{1,\dots ,k\}} = \sum_{I_1 \amalg \dots \amalg I_{\mu} = \{1,\dots ,k \}} (-1)^{\mu -1} (\mu -1)!  S_{I_1}^{\langle i \rangle } \dots S_{I_{\mu }}^{\langle i \rangle },  \nonumber \\
& T_{\{1,\dots ,k\}} =-(-l(l+1))^{k-1}\sum_{i=0}^3 \frac{l_i(l_i+1)}{2} T^{o,\langle i \rangle } _{\{1,\dots ,k\}},  \nonumber
\end{align}
where the sum $\sum_{I_1 \amalg \dots \amalg I_{\mu} = \{1,\dots ,k \}}$ runs over all different partitions of $\{1,\dots ,k\}$. For example
\begin{align}
& T^o_{\phi}=1, \;  T^o_{\{1\}}= S_{\{1\}}, \; T^o_{\{1,2\}}= S_{\{1,2\}}-S_{\{1\}}S_{\{2\}}, \nonumber \\
& T^o_{\{1,2,3\}}= S_{\{1,2,3\}}-S_{\{1\}}S_{\{2,3\}}-S_{\{2\}}S_{\{1,3\}}-S_{\{3\}}S_{\{1,2\}} +2 S_{\{1\}}S_{\{2\}}S_{\{3\}}.\nonumber 
\end{align}
Set
\begin{align}
& \Delta _{\{1,\dots ,k\}} = \sum_{0\leq j \leq [\frac{k}{2}]} \frac{(l(l+1))^j}{2^kj!(k-2j)!} \sum_{w\in W(B_k)} \epsilon (w)w\left( \wp (x_1-x_2)\wp (x_3-x_4)\dots \right. \\
& \; \; \; \; \; \; \; \left. \cdot \wp (x_{2j-1}-x_{2j})\frac{\partial }{\partial x_{2j+1}}\frac{\partial }{\partial x_{2j+2}}\dots \frac{\partial }{\partial x_{k}}\right), \nonumber \\
& q _{\{1,\dots ,k\}} = \sum_{I_1 \amalg \dots \amalg I_{\mu} \in \{1,\dots ,k\} } T_{I_1} \dots T_{I_{\mu }}, \nonumber \\
& T_{w(\{1,\dots ,k\})}=w(T_{\{1,\dots ,k\}} ), \; \; \Delta_{w(\{1,\dots ,k\})}=w(\Delta _{\{1,\dots ,k\}} ) ,\; \; {\mbox{for }}w \in {\mathfrak S}_N , \nonumber
\end{align}
where $[\frac{k}{2}]$ represents the maximum integer not greater than $\frac{k}{2}$.

\begin{prop} $($\cite{O} theorem 7.2, remark 7.4$)$.
Set
\begin{align}
& P_{N-k}= \sum_{i=k}^N \sum_{j=i}^N \frac{1}{i!(j-i)!(N-j)!}\sum_{w\in  {\mathfrak S}_N }\sum_{I_1 \amalg \dots \amalg I_{k} = \{1,\dots ,i \}} \label{eq:PNk}\\
& w\left((-l(l+1))^{i-k}2^{-k}T^o_{I_1}\dots T^o_{I_k} q_{\{i+1, \dots ,j\}} \Delta ^2 _{\{j+1,\dots ,N\}}\right) \nonumber
\end{align}
for $k=0, \dots ,N-1$. Then the operators $P_{j}$ $(1 \leq j \leq N)$ are $W(B_N)$-invariant and
\begin{equation}
[P_j, P_k]=0 \label{Pjkcomm}
\end{equation}
for $1 \leq j,k \leq N$. The Hamiltonian $H$ (\ref{InoHam}) is a linear combination of $P_1$ and $1$, i.e. $H=AP_1 +B$ for some constants $A,B$.
\end{prop}

\subsection{Commuting operators and invariant subspaces} \label{sec:commQES} $ $

We change variables by $z_j=\wp (x_j)$ $(j=1, \dots ,N)$ and set
\begin{equation}
\widehat{P}_k= \Phi(z)^{-1} \circ P_k \circ  \Phi(z), \; \; \; (k=1, \dots ,N)
\label{def:widePk}
\end{equation}
where $\Phi(z)$ is defined in (\ref{def:Phi}). From (\ref{Pjkcomm}), we obtain
\begin{equation}
[\widehat{P}_j, \widehat{P}_k]=0 \label{Phatjkcomm}
\end{equation}
for $1 \leq j,k \leq N$. 

\begin{prop} \label{prop:Pkexpr}
The operators $\widehat{P}_k$ $(k=1, \dots ,N)$ admit the following expansion:
\begin{equation}
\widehat{P}_k= \sum_{0\leq i_1, \dots ,i_N \leq 2 \atop{i_1+ \dots +i_N \leq 2k}} c_{i_1, \dots , i_N} (z)\left(\frac{\partial }{\partial z_{1}}\right)^{i_1}  \dots \left(\frac{\partial }{\partial z_{N}}\right)^{i_N}.
\label{expr:wPk}
\end{equation}
Here, the operators $\widehat{P}_k$ are symmetric in $z_1, \dots ,z_N$ and the coefficients $c_{i_1, \dots , i_N} (z)$ are rational functions in $z_1, \dots ,z_N$ which may have poles only along $z_j-e_i=0$ and $z_{j_1}-z_{j_2}=0$ $(1 \leq i \leq 3, \; 1 \leq j, j_1, j_2 \leq N, \; j_1 \neq j_2)$.
\end{prop}
\begin{proof}
From (\ref{eq:PNk}), the operator $P_k$ admits the expansion 
\begin{equation}
P_k= \sum_{0\leq i_1, \dots ,i_N \leq 2 \atop{i_1+ \dots +i_N \leq 2k}} d_{i_1, \dots , i_N}(x_1, \dots ,x_N) \left(\frac{\partial }{\partial x_{1}}\right)^{i_1} \dots \left(\frac{\partial }{\partial x_{N}}\right)^{i_N}
\end{equation}
such that $d_{i_1, \dots , i_N}(x_1, \dots ,x_N) $ is doubly periodic for each variable $x_1, \dots ,x_N$. From the $W(B_N)$-symmetry of the operator $P_k$,
we obtain $d_{i_1,  \dots , i_N} (x_1, \dots ,-x_j, \dots ,x_N) = (-1)^{i_j} d_{i_1,  \dots , i_N} (x_1, \dots ,x_j, \dots ,x_N) $ for $j=1,\dots ,N$.

Since $\wp'(x)$ is an odd doubly periodic function, the function $d^e_{i_1,  \dots , i_N}(x_1,  \dots ,x_N)$ defined by $d_{i_1,  \dots , i_N}(x_1, \dots ,x_N) = d^e_{i_1,  \dots , i_N}(x_1,  \dots ,x_N) \prod_{j  ; \: i_j=1}\wp'(x_j)$ is even doubly periodic in each $x_j$ $(j=1,\dots ,N)$.

Now, change variables $z_j=\wp(x_j)$ $(j=1,\dots ,N)$ and write 
\begin{equation}
P_k= \sum_{0\leq i_1, \dots ,i_N \leq 2 \atop{i_1+ \dots +i_N \leq 2k}} \tilde{d}_{i_1, \dots , i_N}(x_1, \dots ,x_N) \left(\frac{\partial }{\partial z_{1}}\right)^{i_1} \dots \left(\frac{\partial }{\partial z_{N}}\right)^{i_N}.
\end{equation}
Since $\wp'(x)^2$ and $\wp''(x)$ are even doubly periodic, the function $\tilde{d}_{i_1, \dots , i_N}(x_1, \dots ,x_N) $ is even doubly periodic in each $x_j$ $(j=1,\dots ,N)$.
It is known that even doubly periodic function in $x$ is expressed as $F(\wp(x))$ with some rational function $F(z)$. 
Hence if we set $z_j=\wp(x_j)$ $(j=1,\dots ,N)$, the functions $\tilde{d}_{i_1, \dots , i_N}(x_1, \dots ,x_N) $ admit the expression $\tilde{d}_{i_1, \dots , i_N}(x_1, \dots ,x_N)=\tilde{c}_{i_1, \dots , i_N}(\wp(x_1), \dots ,\wp(x_N))  $ with some rational function $\tilde{c}_{i_1, \dots , i_N}(z_1, \dots ,z_N)$.
From equality (\ref{eq:PNk}) and formulae (\ref{wpformula}), it is shown that the coefficients $\tilde{c}_{i_1, \dots , i_N} (z_1, \dots ,z_N)$ may have poles only along $z_j-e_i=0$ and $z_{j_1}-z_{j_2}=0$ $(1 \leq i \leq 3, \; 1 \leq j, j_1, j_2 \leq N, \; j_1 \neq j_2)$. 

Next we consider the coefficients of the operators $\widehat{P}_k$ defined by (\ref{def:widePk}).
From properties of the functions $\tilde{c}_{i_1, \dots , i_N}(z_1, \dots ,z_N) $, it is shown that the coefficients $c_{i_1, \dots , i_N} (z)$ in (\ref{expr:wPk}) are rational functions in variables $z_1, \dots ,z_N$ which may have poles only along $z_j-e_i=0$ and $z_{j_1}-z_{j_2}=0$ $(1 \leq i \leq 3, \; 1 \leq j, j_1, j_2 \leq N, \; j_1 \neq j_2)$. 

Since the operators $P_k$ and the function $\Phi(z)$ are symmetric in $z_1, \dots ,z_N$, the operators $\widehat{P}_k$ are symmetric in $z_1, \dots ,z_N$.
\end{proof}

\begin{thm} \label{thm:Pinv} 
(i) Let $P^{\mbox{\rm \scriptsize sym}}$ is the space of symmetric polynomials in variables $z_1, \dots ,z_N$, then $\widehat{P}_k  \cdot P^{\mbox{\rm \scriptsize sym}} \subset P^{\mbox{\rm \scriptsize sym}}$ for $k=1,2,\dots ,N$\\
(ii) Let $a$, $b_i$ $(i=0,1,2,3)$ be numbers which satisfy $a\in \{ -l, l+1 \}$ and $b_i \in \{-\frac{l_i}{2}, \frac{l_i+1}{2}\}$ $(i=0,1,2,3)$. Assume that $d=-((N-1)a+b_0+b_1+b_2+b_3)$ is a non-negative integer. Then $\widehat{P}_k  \cdot V_d^{\mbox{\rm \scriptsize sym}} \subset V_d^{\mbox{\rm \scriptsize sym}}$ for $k=1,2,\dots ,N$, where $V_d^{\mbox{\rm \scriptsize sym}}$ is the finite-dimensional space defined in proposition \ref{prop:Hinv}.
\end{thm}
\begin{proof}
The coefficients $c_{i_1, \dots , i_N} (z)$ in (\ref{expr:wPk}) are rational functions which may have poles only along $z_j-e_i=0$, $z_{j_1}-z_{j_2}=0$ $(1 \leq i \leq 3, \; 1 \leq j, j_1, j_2 \leq N, \; j_1 \neq j_2)$.

Let us fix $j_1$ and $j_2$ that satisfy $j_1\neq j_2$. Let $p$ be the maximal number of degrees of a pole along $z_{j_1}-z_{j_2}=0$ of functions $c_{i_1, \dots , i_N} (z)$ for possible all $i_1, \dots , i_N$.
If $f (z) \in  P^{\mbox{\scriptsize sym}}$ then the function $\widehat{P}_k f(z)$ has a pole at most degree $p$ along $z_{j_1}-z_{j_2}=0$.
Let $f(z)$ be an element of $P^{\mbox{\scriptsize sym}}$ such that the function $\widehat{P}_k f(z)$ has a pole of maximum degree along $z_{j_1}-z_{j_2}=0$. We denote the degree by $p'$. Then the function $\widehat{P}_k \widehat{H} f(z)$ has a pole of degree at most $p'$ along $z_{j_1}-z_{j_2}=0$, because $\widehat{H} f(z) \in P^{\mbox{\scriptsize sym}}$. On the other hand, it is shown that the function $\widehat{H} \widehat{P}_k f(z)$ has a pole of degree $p'+2$ along $z_{j_1}-z_{j_2}=0$ if $p' \neq 0, \: 1-2a$.
Since $\widehat{P}_k \widehat{H}f(z)= \widehat{H} \widehat{P}_k f(z)$, if $1-2a \not \in \Zint$ then $p'$ must be equal to zero. Hence the function $\widehat{P}_k f(z)$ is holomorphic along $z_{j_1}-z_{j_2}=0$ if $1-2a \not \in \Zint$.
By a continuity argument in $a$, it is shown that the function $\widehat{P}_k f(z)$ is holomorphic along $z_{j_1}-z_{j_2}=0$ for all $a$.

Similarly if $f(z) \in  P^{\mbox{\scriptsize sym}}$ then the function $\widehat{P}_k f(z)$ is holomorphic along $z_j-e_i=0$ and $z_{j_1}-z_{j_2}=0$ $(1 \leq i \leq 3, \; 1 \leq j, j_1, j_2 \leq N, \; j_1 \neq j_2)$. Hence $\widehat{P}_k f(z) \in  P^{\mbox{\scriptsize sym}}$.

Next we prove (ii).
%

From expression (\ref{expr:wPk}), there exists $p\in \Zint_{\geq 0}$ such that $\widehat{P}_k f(z) \subset  V_{d+p}^{\mbox{\scriptsize sym}}$ for all $f(z)\in  V_d^{\mbox{\scriptsize sym}}$. Let $f(z)$ be an element of $V_d^{\mbox{\scriptsize sym}}$ such that the degree of $\widehat{P}_k f(z)$ is maximum.
We denote the degree by $d+p'$.
Then the function $\widehat{P}_k \widehat{H} f(z)$ has a degree at most $d+p'$, because $\widehat{H} f(z) \in  V_d^{\mbox{\scriptsize sym}}$. On the other hand it is shown that the function $\widehat{H} \widehat{P}_k f(z)$ has a degree $d+p'+1$ if $p' \neq 0, \: 2b_0-1/2$.
From the commutativity of $\widehat{H}$ and $\widehat{P}_k$, if $ 2b_0-1/2 \not \in \Zint$ then $p'$ must be equal to zero. Thus $\widehat{P}_k f(z) \in V_d^{\mbox{\scriptsize sym}}$. By a continuity argument, we can remove the condition $ 2b_0-1/2 \not \in \Zint$.

Hence we obtain(ii).
\end{proof}

In summary, we established that the higher commuting Hamiltonians also preserve the space related to the quasi-exact solvablity in theorem \ref{thm:Pinv} (ii).

\subsection{Relationship to the $L^2$ space} \label{sec:L2} $ $

Assume $l,l_0,l_1 \in \Rea _{\geq 0}$ and $l_2, l_3 \in \Rea $ in this subsection.

The invariant space $V_d^{\mbox{\scriptsize sym}}$ of the $BC_N$ Inozemtsev model is defined for each $a\in \{ -l, l+1 \}$ and $b_i \in \{-\frac{l_i}{2}, \frac{l_i+1}{2}\}$ $(i=0,1,2,3)$ with the condition $d=-((N-1)a+b_0+b_1+b_2+b_3) \in \Zint_{\geq 0}$. For these numbers $a, b_0, b_1, b_2, b_3$, define
\begin{equation}
W_d^{\mbox{\scriptsize sym}}=\{ \Phi(\wp(x_1), \dots , \wp(x_N)) f(\wp(x_1), \dots , \wp(x_N)) | f(z_1, \dots ,z_N) \in V_d^{\mbox{\scriptsize sym}} \}, \label{Wdsym}
\end{equation}
where the function $\Phi(z_1, \dots , z_N) =\prod_{1\leq j<k\leq N} (z_j-z_k)^a \prod_{j=1}^N \prod _{i=1}^3 (z_j-e_i)^{b_i}$ was defined in (\ref{def:Phi}).

From relations (\ref{def:Phi}), (\ref{def:widePk}), and Theorem \ref{thm:Pinv}, the following proposition is shown immediately:

\begin{prop} \label{prop:PinvW}
Assume that $d=-((N-1)a+b_0+b_1+b_2+b_3)$ is a non-negative integer. Then $H  \cdot W_d^{\mbox{\rm \scriptsize sym}} \subset W_d^{{\mbox{\rm \scriptsize sym}}}$ and $P_k  \cdot W_d^{{\mbox{\rm \scriptsize sym}}} \subset W_d^{{\mbox{\rm \scriptsize sym}}}$ for $k=1,2,\dots ,N$.
\end{prop}

We look into the condition that the space $W_d^{{\mbox{\scriptsize sym}}}$ lies in $L^2$ space.
If $f(x) \in W_d^{{\mbox{\scriptsize sym}}}$, then $|f(x)| \sim |x_j-x_k|^{a}$ (resp. $|f(x)| \sim |x_j+x_k|^{a}$)  as $|x_j-x_k| \rightarrow 0$ (resp. $|x_j+x_k| \rightarrow 0$) $(j\neq k)$, $|f(x)| \sim |x_j|^{2b_0}$ as $|x_j| \rightarrow 0$, and $|f(x)| \sim |x_j-\frac{1}{2}|^{2b_1}$ as $|x_j-\frac{1}{2}| \rightarrow 0$. 
Hence if $a=l+1$, $b_0 =\frac{l_0+1}{2}$, and $b_1= \frac{l_1+1}{2}$ then the function $f(x)\in W_d^{{\mbox{\scriptsize sym}}}$ is locally square-integrable on $|x_j-x_k| , |x_j+x_k| , |x_j|, |x_j-\frac{1}{2}| < \epsilon$ for sufficiently small $\epsilon $.
Combining with the periodicity of $f(x)\in W_d^{{\mbox{\scriptsize sym}}}$, we obtain $\int_{0< x_1< \dots < x_N < 1}|f(x)|^2dx_1\dots dx_N < \infty$. Hence the following proposition is shown:

\begin{prop}
Let $b_i \in \{-\frac{l_i}{2}, \frac{l_i+1}{2}\}$ $(i=2,3)$. If $d= -((N-1)(l+1)+\frac{l_0+l_1}{2}+1+b_2+b_3) \in \Zint_{\geq 0}$ then every function in $W_d^{{\mbox{\rm \scriptsize sym}}}$ is square-integrable on the domain  $0< x_1< \dots <x_N < 1$.
\end{prop}

In the case $d= -((N-1)(l+1)+\frac{l_0+l_1}{2}+1+b_2+b_3) \in \Zint_{\geq 0}$, some eigenvalues of the commuting Hamiltonians on the Hilbert space ($L^2$ space) appear as the eigenvalues on the subspace $W_d^{{\mbox{\scriptsize sym}}}$. Hence some eigenvalues on the Hilbert space would be obtained explicitly, because the eigenvalues in the finite-dimensional space are got by algebraic calculations. 
Note that the case $BC_1$ was done in \cite{Tak2}, and Gomez-Ullate, Gonzalez-Lopez, and Rodriguez considered the relationship between the $L^2$ space and the space related to the quasi-exact solvablity for some special cases in \cite{GGR}.

As an aside, the joint eigenvalues of the trigonometric $BC_N$ Calogero-Sutherland model are already known and their expression is simple.
Distributions of eigenvalues in $L^2 \cap W_d^{{\mbox{\scriptsize sym}}}$ will be detected by considering the trigonometric limit $\tau \rightarrow \sqrt{-1} \infty$ while fixing coupling constants $l, l_0, l_1, l_2, l_3$.

\subsection{The case $l, l_0, l_1, l_2, l_3 \in \Zint_{\geq 0}$} \label{sec:int} $ $

Let us consider the case $l, l_0, l_1, l_2, l_3 \in \Zint_{\geq 0}$. In this case, the Hamiltonian and the higher commuting Hamiltonians preserve several finite-dimensional spaces of elliptic functions.

The invariant space $W_d^{{\mbox{\scriptsize sym}}}$ is defined for each $a\in \{ -l, l+1 \}$ and $b_i \in \{-\frac{l_i}{2}, \frac{l_i+1}{2}\}$ $(i=0,1,2,3)$ with the condition $d=-((N-1)a+b_0+b_1+b_2+b_3) \in \Zint_{\geq 0}$. 
 
For each $l,  l_0, l_1, l_2, l_3 \in \Zint_{\geq 0}$, there are eight possible sets of $(a, b_0, b_1, b_2, b_3)$ for each set  the invariant space $W_d^{{\mbox{\scriptsize sym}}}$ is defined, if $N\geq 2$. If $N=1$ then there are four possible sets of $( b_0, b_1, b_2, b_3)$.
For example, if $N\geq 2$, $l\gg \l_0, l_1, l_2, l_3$ and $(N-1)l+l_0+l_1+l_2+l_3 \in 2\Zint _{>0}$, then the cases
$$
(a, b_0, b_1, b_2, b_3)= \left\{
\begin{array}{l}
(-l, \frac{-l_0}{2},  \frac{-l_1}{2},  \frac{-l_2}{2},  \frac{-l_3}{2}), \: (-l, \frac{l_0+1}{2},  \frac{l_1+1}{2},  \frac{l_2+1}{2},  \frac{l_3+1}{2}), \\
(-l, \frac{-l_0}{2},  \frac{-l_1}{2},  \frac{l_2+1}{2},  \frac{l_3+1}{2}), \: (-l, \frac{l_0+1}{2},  \frac{l_1+1}{2},  \frac{-l_2}{2},  \frac{-l_3}{2}), \\
(-l, \frac{-l_0}{2},  \frac{l_1+1}{2},  \frac{-l_2}{2},  \frac{l_3+1}{2}), \: (-l, \frac{l_0+1}{2},  \frac{-l_1}{2},  \frac{l_2+1}{2},  \frac{-l_3}{2}), \\
(-l, \frac{-l_0}{2},  \frac{l_1+1}{2},  \frac{l_2+1}{2},  \frac{-l_3}{2}), \: (-l, \frac{l_0+1}{2},  \frac{-l_1}{2},  \frac{-l_2}{2},  \frac{l_3+1}{2}), 
\end{array} \right.
$$
are permitted.

By a straightforward calculation, the dimension of direct sum of spaces $W_d^{{\mbox{\scriptsize sym}}}$ of elliptic functions can be calculated. For the case $N=1$, the dimension is
$$
\left\{ 
\begin{array}{ll}
2k_0+1, & \tilde{l} \mbox{ is even and } k_0+k_3\geq \frac{\tilde{l}}{2};\\
\tilde{l}-2k_3+1, & \tilde{l} \mbox{ is even and } k_0+k_3< \frac{ \tilde{l}}{2};\\
2k_0+1, & \tilde{l} \mbox{ is odd and } k_0\geq \frac{ \tilde{l} +1}{2};\\
\tilde{l}+2, & \tilde{l} \mbox{ is odd and } k_0< \frac{ \tilde{l} +1}{2},
\end{array}
\right.
$$
where $\tilde{l}= l_0+l_1+l_2+l_3$, $k_0= \max (l_0, l_1, l_2, l_3)$, and $k_3=\min (l_0, l_1, l_2, l_3)$ (see also \cite{Tak1,Tak3}). For the case $N=2$, the dimension is $(2l+1)^2+\sum_{i=0}^3 l_i(l_i+1)$ for all the case $l,  l_0, l_1, l_2, l_3 \in \Zint_{\geq 0}$. Since the dimension is directly related to the genus of the spectral curve for the case $N=1$ \cite{Tak3}, the dimension for the case $N\geq 2$ might also play important roles.

\section{Ruijsenaars models and Inozemtsev models} \label{sec:Ruij}

In \cite{Rui}, Ruijsenaars introduced a relativistic version of the Calogero-Moser-Sutherland model, which is called the Ruijsenaars model of type $A_N$ these days.
In \cite{vD}, van Diejen introduced the $BC_N$ Ruijsenaars-type model which has ten parameters $(\kappa; \mu, \nu_0, \bar{\nu}_0, \nu_1, \bar{\nu}_1, \nu_2, \bar{\nu}_2, \nu_3, \bar{\nu}_3)$, and Hikami and Komori constructed higher commuting operators by use of root algebra in \cite{HK,Kom1,Kom2,KH}, that ensures the integrability.

The lowest operator of the $BC_N$ (or $A^{(2)}_{2N}$) Ruijsenaars model is given as follows: 
\begin{align}
Y_1=& \sum_{j=1}^N\left( \prod_{k=1\atop{k\neq j}}^N\frac{\theta _1(x_j-x_k-\mu)}{\theta _1(x_j-x_k)}\frac{\theta _1(x_j+x_k-\mu)}{\theta _1(x_j+x_k)}\right) \\
& \quad \cdot \left( \prod_{r=0}^3 \frac{\theta_{r+1}(x_j-\nu_r)}{\theta_{r+1}(x_j)}\frac{\theta_{r+1}(x_j+\kappa/2-\bar{\nu}_r)}{\theta_{r+1}(x_j+\kappa/2)}\right) t_j(\kappa) \nonumber \\
& + \sum_{j=1}^N\left( \prod_{k=1\atop{k\neq j}}^N\frac{\theta _1(x_j+x_k+\mu)}{\theta _1(x_j+x_k)}\frac{\theta _1(x_j-x_k+\mu)}{\theta _1(x_j-x_k)}\right) \nonumber \\
& \quad \cdot \left( \prod_{r=0}^3 \frac{\theta_{r+1}(x_j+\nu_r)}{\theta_{r+1}(x_j)}\frac{\theta_{r+1}(x_j-\kappa/2+\bar{\nu}_r)}{\theta_{r+1}(x_j-\kappa/2)}\right) t_j(-\kappa) \nonumber \\
& + \sum_{p=0}^3 \left( \frac{\pi}{\theta ' _1 (0)} \right)^2\frac{2}{\theta_1(\mu)\theta_1(\kappa+\mu)}\left(\prod_{r=0}^3\theta_{r+1}(\kappa/2+\nu_{\pi_p r})\theta_{r+1}(\bar{\nu}_{\pi_p r})\right) \nonumber \\
& \quad \cdot \left( \prod_{j=1}^N \frac{\theta_{p+1}(x_j-\kappa /2 -\mu )}{\theta_{p+1}(x_j-\kappa /2 )} \frac{\theta_{p+1}(x_j+\kappa /2 +\mu )}{\theta_{p+1}(x_j+\kappa /2 )} \right). \nonumber
\end{align}
Here $\theta_j(x)$ $(j=1,2,3,4)$ is the Jacobi theta function (see (\ref{def:theta})) and $t_i(\kappa)$ is a translation of the variable $x_i$ by $\kappa$, i.e. $t_i(\kappa) f(x_1, \dots ,x_i, \dots x_N)= f(x_1, \dots ,x_i+\kappa , \dots x_N)$. $\pi _r$ $(r=0,1,2,3)$ denotes the permutation $\pi_0=id$, $\pi_1=(01)(23)$, $\pi_2=(02)(13)$, $\pi_3=(03)(12)$, where $(ij)k= \left\{ \begin{array}{ll}
k, & k\neq i,j; \\
j, & k=i; \\
i, & k=j. 
\end{array}
\right. $


Hikami and Komori showed that if $k= \left( 2(N-1)\mu +\sum_{i=0}^3(\nu_i +\bar{\nu }_i) \right)/\kappa \in 2\Zint_{\geq 0}$ then the operator $Y_1$ and higher commuting operators preserve the space of level $k$ theta functions of type $A^{(2)}_{2N}$. They proved it using root algebra. Their presentation of the invariant subspace would be technical for non-experts. In this section we describe them plainly.

The space of level $k$ theta functions is defined as follows:
\begin{equation}
Th_k^{W(B_N)}= \left\{ f: \Cplx ^N \rightarrow \Cplx \left| 
\begin{array}{l}
\mbox{ holomorphic, }W(B_N)\mbox{-invariant} \\
f(x+n )=f(x) , \quad \quad (\forall n \in \Zint^N ) \\
f(x+n\tau)=f(x) e^{-2\pi \sqrt{-1} k((x|n)+(n|n)\tau/2)}
\end{array}
\right. \right\} , \label{ThkW}
\end{equation}
where $(x|y)=\sum_{i=1}^Nx_iy_i $ for $x=(x_1,\dots ,x_N) \in \Cplx ^N$ and  $y=(y_1,\dots ,y_N) \in \Cplx ^N$.
A function $f(x_1, \dots ,x_N)$ is $W(B_N)$-invariant if and only if the relations $f(x_{\sigma (1)}, \dots ,x_{\sigma (N)})=f(x_1, \dots ,x_N)$ for $\forall \sigma \in {\frak S}_N$ and $f(x_1, \dots ,x_i, \dots ,x_N)=f(x_1, \dots ,-x_i, \dots ,x_N)$ for $\forall i \in \{1, \dots ,N \}$ are satisfied. For the case $N=1$, we obtain  $\dim Th_{2l}^{W(B_1)} =l+1$ for $l \in \Zint_{\geq 0}$. Let $\theta ^{(1)}(x), \dots ,\theta ^{(l+1)}(x)$ be a basis of $Th_{2l}^{W(B_1)}$.
Then the space $Th_{2l}^{W(B_N)}$ is spanned by functions $\sum _{\sigma \in {\frak S}_N} \theta ^{(k_1)}(x_{\sigma(1)}) \dots \theta ^{(k_N)}(x_{\sigma(N)})$ $(1\leq k_1 \leq \dots \leq k_N \leq l+1)$.
Therefore, $\dim Th_{2l}^{W(B_N)}
=\frac{(l+N)!}{l!N!}$ for $l \in \Zint_{\geq 0}$.

\begin{prop} $($c.f. \cite{Kom2}$)$
If $k= \left( 2(N-1)\mu +\sum_{i=0}^3(\nu_i +\bar{\nu }_i) \right)/\kappa \in 2\Zint_{\geq 0}$ then the operator $Y_1$ preserve the space $Th_k^{W(B_N)}$.
\end{prop}
\begin{proof}
Let $f(x) \in Th_k^{W(B_N)}$. Then the function $Y_1 f(x)$ is $W(B_N)$-invariant. From the quasi-periodicity of $\theta _i(x)$ $(i=0,1,2,3)$ (see (\ref{thetaperiod})) and $f(x)$ (see (\ref{ThkW})), the function $Y_1 f(x)$ has a quasi-periodicity as condition (\ref{ThkW}) when $ k= ( 2(N-1)\mu +\sum_{i=0}^3(\nu_i +\bar{\nu }_i))/\kappa $. Hence if we show the holomorphy of the function $Y_1 f(x)$ on $\Cplx ^N$, we have $Y_1 f(x) \in Th_k^{W(B_N)}$.
Thus it is sufficient to show that the residues of the function $Y_1 f(x)$ at $x_j-x_k=0$, $x_j+x_k=0$ $(1\leq j\neq  k \leq N)$, and $x_j=0, 1/2, (1+\tau)/2, \tau/2 , \pm \kappa /2, 1/2\pm \kappa/2 , (1+\tau)/2\pm \kappa /2, \tau/2\pm \kappa /2$ $(1\leq j \leq N)$ are zero. These are shown directly by using the quasi-periodicity of $f(x)$ (\ref{ThkW}) and  $\theta _i(x)$ $(i=0,1,2,3)$ (\ref{thetaperiod}). Note that we rely on the condition $k \in 2\Zint $ in this step.
\end{proof}

Let us consider the non-relativistic (difference-differential) limit of the Ruijsenaars model. It is known that the Inozemtsev model appear by this limit. Now we will exhibit it explicitly.

Let $a=-\mu /\kappa $, $b_0=-(\nu_1 +\bar{\nu }_1) /2\kappa $, $b_1=-(\nu_2 +\bar{\nu }_2 )/2\kappa $, $b_2=-(\nu_3 +\bar{\nu }_3 )/2\kappa $, $b_3=-(\nu_0 +\bar{\nu }_0 )/2\kappa $ and
\begin{equation}
\Theta (x)=\prod_{1\leq j<k\leq N} \left( \theta_1(x_j-x_k)\theta_1(x_j+x_k) \right)^{a}  \prod _{j=1}^N \theta _1(x_j)^{2b_0}\theta _2(x_j)^{2b_1}\theta _3(x_j)^{2b_2}\theta _0(x_j)^{2b_3}.
\end{equation}
Assume $a \in \{-l, l+1 \}$, $b_0 \in \{-l_0/2, (l_0+1)/2\}$, $b_1 \in \{-l_1/2, (l_1+1)/2\}$, $b_2 \in \{-l_2/2, (l_2+1)/2\}$, and $b_3 \in \{-l_3/2, (l_3+1)/2\}$.
 
As $\kappa \rightarrow 0$ while $a$, $b_0$, $b_1$, $b_2$, $b_3$ are fixed,
\begin{equation}
( -\Theta (x)\circ Y_1 \circ \Theta (x) ^{-1} +C_0)/\kappa ^2 \rightarrow H, \label{Thetax}
\end{equation}
where $H$ is the Hamiltonian of the $BC_N$ Inozemtsev model given in (\ref{InoHam}) and $C_0$ is a constant. Hence, we recover the Hamiltonian of the $BC_N$ Inozemtsev model from a operator of $BC_N$ Ruijsenaars model via a limit $\kappa \rightarrow 0$.

Let us make a correspondence between the invariant spaces of theta functions on the $BC_N$ Ruijsenaars model and the space related to the quasi-exact solvability on the $BC_N$ Inozemtsev model.

\begin{prop} \label{prop:isom}
Let   $a=-\mu /\kappa $, $b_0=-(\nu_1 +\bar{\nu }_1) /2\kappa $, $b_1=-(\nu_2 +\bar{\nu }_2 )/2\kappa $, $b_2=-(\nu_3 +\bar{\nu }_3 )/2\kappa $, $b_3=-(\nu_0 +\bar{\nu }_0 )/2\kappa $.
Let $Th^{W(B_N)}_{2k}$ (\ref{ThkW}) be the theta-type invariant space of $BC_N$ Ruijsenaars model,  $W^{{\mbox{\rm \scriptsize sym}}}_k$ (\ref{Wdsym}) be the invariant space of $BC_N$ Inozemtsev model, and $\Theta (x) $ be the function defined in (\ref{Thetax}).
Assume $k= - ((N-1)a +b_0+b_1+b_2+b_3) \in \Zint_{\geq 0}$.

Then the map 
\begin{equation}
\begin{array}{cccc}
\phi  : & Th^{W(B_N)}_{2k} & \rightarrow & W^{{\mbox{\rm \scriptsize sym}}}_k \\
          & f(x_1,\dots x_N) & \mapsto &  \Theta (x)  f(x_1,\dots x_N)
\end{array}
\end{equation}
is an isomorphism of vector spaces.
\end{prop}

\begin{proof}
Let us consider the correspondence between the space $Th^{W(B_N)}_{2k}$ and the space $V^{{\mbox{\scriptsize sym}}}_k$, where $V^{\mbox{\scriptsize sym}}_k$ was defined in proposition \ref{prop:Hinv}.

Let $f(x_1, \dots ,x_N) \in Th^{W(B_N)}_{2k}$ and $ g(x_1, \dots ,x_N) $ $\! = f (x_1, \dots ,x_N) $ $\! \Theta (x) \Phi(\wp(x_1), \dots , \wp(x_N))^{-1}$, where the function $\Phi(z_1, \dots , z_N)$ was defined in (\ref{def:Phi}).
From the condition $k \in \Zint$, the function $g(x_1, \dots ,x_N)$ does not have branchs on $\Cplx ^N$. It is seen that the function $g(x_1, \dots ,x_N)$ is doubly periodic, $W(B_N)$-invariant, and may have poles only along $x_j=0$ $(j=1,\dots ,N)$ up to periods with degree at most $k$.

Hence there exists $ \tilde{g}(z_1, \dots ,z_N) \in V^{{\mbox{\scriptsize sym}}}_k$ such that $\tilde{g}(\wp(x_1), \dots ,\wp(x_N)) =g(x_1, \dots ,x_N) $ by a similar argument in the proof of proposition \ref{prop:Pkexpr}.

By composing with the canonical map from $V^{{\mbox{\scriptsize sym}}}_k$ to $W^{{\mbox{\scriptsize sym}}}_k$, we obtain $\phi (Th^{W(B_N)}_{2k}) \subset W^{{\mbox{\scriptsize sym}}}_k$. It is obvious that the map $\phi$ is injective, and the dimension of $Th^{W(B_N)}_{2k}$ is equal to that of $W^{{\mbox{\scriptsize sym}}}_k$. Therefore the map $\phi$ is bijective.
\end{proof}

In proposition \ref{prop:isom}, we have established that the theta-type invariant space $Th^{W(B_N)}_{2k}$ of the $BC_N$ Ruijsenaars model corresponds to the space $W^{{\mbox{\scriptsize sym}}}_k$ which is related to the quasi-exact solvability of the $BC_N$ Inozemtsev model.

\section{Degenerate Inozemtsev model} \label{sec:deg}

\subsection{Trigonometric $BC_N$ Inozemtsev model} $ $

In \cite{ST}, Sasaki and Takasaki considered degenerate $BC_N$ Inozemtsev models and their quasi-exact solvability. They also considered for the case of type $A_N$.

In this section, we consider the degeneration of $BC_N$ Inozemtsev model and show that the finite-dimensional invariant spaces for the elliptic $BC_N$ Inozemtsev model tend to the spaces introduced by Sasaki and Takasaki by the degeneration.

The Hamiltonian of the trigonometric (or degenerate) $BC_N$ Inozemtsev model is given as follows:
\begin{align}
& H^{\mbox{\scriptsize (D)}}=-\sum_{j=1}^N\frac{\partial ^2}{\partial x_j^2}+2 l(l+1)\sum_{1\leq j<k\leq N} \left( \frac{\pi ^2}{\sin ^2\pi (x_j-x_k)} +\frac{\pi ^2}{\sin ^2 \pi (x_j +x_k)} \right) \label{triInoHam} \\
& \; \; \;+ \sum_{j=1}^N \left( \frac{\pi^2 l_0(l_0+1)}{\sin ^2 \pi x_j}+ \frac{\pi^2 l_1(l_1+1)}{\cos ^2 \pi x_j}+\tilde{c}_1 \cos 2\pi x_j +\tilde{c}_2\cos 4\pi x_j \right), \nonumber
\end{align}
where $l$, $l_0$, $l_1$, $\tilde{c}_1$, $\tilde{c}_2$ are coupling constants.

This model is known to be quantum integrable. In \cite{O}, Oshima gave the explicit expression of commuting operators of conserved quantities.

Set
\begin{align}
 \Phi_D(x)=& \left| \exp\left( -\frac{\tilde{a}}{2}\sum_{j=1}^N \cos 2\pi x_j \right) \prod_{j=1}^N (\sin \pi x_j)^{l_0+1} (\cos \pi x_j)^{l_1+1}\right. \\
&  \left. \prod_{1\leq j_1<j_2\leq N} \left( \sin \pi (x_{j_1}-x_{j_2}) \sin \pi (x_{j_1 } +x_{j_2}) \right) ^{l+1} \right| \nonumber .
\end{align}
Let $W^{\mbox{\scriptsize (D)}}_L$ be the vector space spanned by functions $\Phi_D(x) \!$ $ (\sin  \pi x_1)^{2m_1} \!$ $(\sin  \pi x_2)^{2m_2} \dots (\sin  \pi x_N)^{2m_N} $ such that $m_i \in \{ 0,1, \dots ,L\}$ for all $i$, and $W^{\mbox{\scriptsize (D),sym}}_L$ be the set of ${\frak S}_N$-invariant elements in $W^{\mbox{\scriptsize (D)}}_L$.

The following proposition is essentially shown in \cite{ST}.
\begin{prop} \label{prop:deginv} $($c.f. \cite{ST}$)$
The Hamiltonian $H^{\mbox{\rm \scriptsize (D)}}$ (\ref{triInoHam}) preserves the space $W^{\mbox{\rm \scriptsize (D),sym}}_L$, if $L\in \Zint _{\geq 0}$, $\tilde{c}_2= -\frac{\pi ^2 \tilde{a}^2}{2}$ and $\tilde{c}_1= 2\tilde{a} \pi^2 ( 2L+l_0+l_1+3+2(N-1)(l+1))$.
\end{prop}
\begin{proof}
We set $W_0= \log \Phi_D (x)$.
Then we have 
\begin{align}
& \sum_{j=1}^N  \left( \left( \frac{\partial W_0}{\partial x_j}\right) ^2+ \frac{\partial ^2 W_0}{\partial x_j^2 } \right) = \sum_{1\leq j<k\leq N} \left( \frac{2\pi ^2 l(l+1)}{\sin ^2\pi (x_j-x_k)} +\frac{2\pi ^2 l(l+1)}{\sin ^2 \pi (x_j +x_k)} \right) \\
& \; \; \;+ \sum_{j=1}^N \left( \frac{\pi^2 l_0(l_0+1)}{\sin ^2 \pi x_j}+ \frac{\pi^2 l_1(l_1+1)}{\cos ^2 \pi x_j}+\tilde{c}_3 \cos 2\pi x_j -\frac{\pi ^2 \tilde{a}^2}{2} \cos 4 \pi x_j \right) +C_0, \nonumber
\end{align}
where $C_0$ is a constant term and $\tilde{c}_3 = 2\tilde{a} \pi^2 ( l_0+l_1+3+2(N-1)(l+1))$. By comparing with the Hamiltonian in \cite[(7.1)]{ST} and its corresponding 'exactly solvable sector' \cite[(7.13)]{ST}, we obtain the proposition.
\end{proof}

\subsection{Degeneration} $ $

In this section, we consider the degeneration (the trigonometric limit) $\tau \rightarrow \sqrt{-1}\infty $ and see the correspondences of Hamiltonians and their invariant spaces between the nondegenerate model and the degenerate one.

Let $l, l_0, l_1, l_2, l_3$ be the coupling constants of elliptic Inozemtsev model (see (\ref{InoHam})). We adopt the following limits of coupling constants as $\tau \rightarrow \sqrt{-1}\infty $:\\
$\bullet$ $l, l_0, l_1$: fixed;\\
$\bullet$ $l_2 = \frac{\tilde{a}}{8} p^{-1}+ \tilde{b}$ and  $l_3 = -\frac{\tilde{a}}{8} p^{-1}+ \tilde{b}$, where $p=\exp (\pi \sqrt{-1}\tau )$.
Here we note that $p \rightarrow 0$ as $\tau \rightarrow \sqrt{-1}\infty $.
Then the Hamiltonian $H$ of the elliptic Inozemtsev model (see (\ref{InoHam})) tends to the Hamiltonian $H^{\mbox{\scriptsize (D)}}$ of the trigonometric Inozemtsev model (see (\ref{triInoHam})). More precisely,
\begin{align}
& H +\frac{\pi^2}{3} \left( l_2(l_2+1) + l_3(l_3+1) \right) \rightarrow -\sum_{j=1}^N\frac{\partial ^2}{\partial x_j^2} \label{limHamil} \\
& \; \; \; \; +\sum_{1\leq j<k\leq N} \left( \frac{2\pi ^2l(l+1)}{\sin ^2\pi (x_j-x_k)} +\frac{2\pi ^2l(l+1)}{\sin ^2 \pi (x_j +x_k)} \right) + \sum_{j=1}^N \left( \frac{\pi^2 l_0(l_0+1)}{\sin ^2 \pi x_j} \right. \nonumber \\
& \: \: \: \: \left. + \frac{\pi^2 l_1(l_1+1)}{\cos ^2 \pi x_j}+2\pi^2 \tilde{a}(2\tilde{b} +1) \cos 2\pi x_j -\frac{\pi ^2 \tilde{a}^2}{2} \cos 4 \pi x_j \right)+ C_1, \nonumber
\end{align}
as $p\rightarrow 0$, where $C_1$ is a constant.

Let us observe how the invariant space varies as $p\rightarrow 0$.
Since $\tilde{b}=\frac{l_2+l_3}{2}$, the Hamiltonian $H$ (see (\ref{InoHam})) preserves the space $W_L^{{\mbox{\scriptsize sym}}}$ (see (\ref{Wdsym})) for $a=l+1, b_1= \frac{l_1+1}{2},  b_2= \frac{-l_2}{2}$, and $b_3= \frac{-l_3}{2} $ if $L=-(N-1)(l+1) -\frac{l_0+l_1+2}{2}+\tilde{b} \in \Zint_{\geq 0}$.

We consider the limit $p \rightarrow 0 $. Note that if $ L=-(N-1)(l+1) -\frac{l_0+l_1+2}{2}+\tilde{b} \in \Zint_{\geq 0}$ then the Hamiltonian $H$ preserves the space $W_L^{{\mbox{\scriptsize sym}}}$ whenever $p$ varies.

Let $\Phi(z)=\prod_{1\leq j<k\leq N} (z_j-z_k)^{l+1} \prod_{j=1}^N  (z_j-e_1)^{\frac{l_1+1}{2}}(z_j-e_2)^{\frac{-l_2}{2}}(z_j-e_3)^{\frac{-l_3}{2}}$ be the function defined in (\ref{def:Phi}) for $a=l+1, b_1= \frac{l_1+1}{2},  b_2= \frac{-l_2}{2}, b_3= \frac{-l_3}{2} $. Then $\Phi( \wp(x_1), \dots , \wp(x_N)) \rightarrow  C_3 \Psi_D(x)$ as $\tau \rightarrow \sqrt{-1}\infty $, where $C_3$ is a constant and
\begin{align}
& \Psi_D(x)=  \prod_{j=1}^N (\sin \pi x_j)^{-2(N-1)(l+1)-(l_1+1)+2\tilde{b} } (\cos \pi x_j)^{l_1+1} \\
& \; \; \; \; \prod_{1\leq j_1<j_2\leq N} \left( \sin \pi (x_{j_1}-x_{j_2}) \sin \pi (x_{j_1 } +x_{j_2}) \right) ^{l+1}  \exp\left( -\frac{\tilde{a}}{2}\sum_{j=1}^N \cos 2\pi x_j \right) \nonumber.
\end{align}

Set $t(x)=\frac{\pi^2}{\sin^2 \pi x}-\frac{\pi ^2}{3}$. Let $\tilde{W}^{\mbox{\scriptsize (D)}}_L$ be the vector space spanned by functions $\Psi_D(x) \!$ $ t(x_1)^{m_1} \!$ $t(x_2)^{m_2} \dots t(x_N)^{m_N} $ such that $m_i \in \{ 0,1, \dots ,L\}$ for all $i$, and $\tilde{W}^{\mbox{\scriptsize (D),sym}}_L$ be the set of ${\frak S}_N$-invariant elements in $\tilde{W}^{\mbox{\scriptsize (D)}}_L$.

As $p\rightarrow 0$, the vector space $W_L^{{\mbox{\scriptsize sym}}}$ tends to the space $\tilde{W}^{\mbox{\scriptsize (D),sym}}_L$, and the operator which appears on the right-hand side of (\ref{limHamil}) preserves the space $\tilde{W}^{\mbox{\scriptsize (D),sym}}_L$ if $L=-(N-1)(l+1) -\frac{l_0+l_1+2}{2}+\tilde{b} \in \Zint_{\geq 0}$.
If $l_0+1=2\tilde{b}-2(N-1)(l+1)-(l_1+1)$ then it is seen that $\tilde{W}^{\mbox{\scriptsize (D),sym}}_L = W^{\mbox{\scriptsize (D),sym}}_L$. 
Therefore, we recover proposition \ref{prop:deginv} by the trigonometric limit.

In summary, by the trigonometric limit we have shown that some finite-dimensional invariant space of the Hamiltonian $H$ of elliptic model (see (\ref{InoHam})) tends to the the invariant space of Sasaki and Takasaki which is related to the quasi-exact solvablity.

In \cite{O}, Oshima described the limit procedure of the commuting operators of conserved quantities. By applying Oshima's result and Proposition \ref{prop:PinvW} in this article, it follows that the commuting operators of the trigonometric $BC_N$ Inozemtsev model also preserve the space $W^{\mbox{\scriptsize (D),sym}}_L$.

Hence we established that the commuting operators of the trigonometric $BC_N$ Inozemtsev model also preserve the space related to the quasi-exact solvability.

\section{Concluding remarks}
In this paper, quasi-exact solvability for the $BC_N$ Inozemtsev model is proved not only for the Hamiltonian but also for commuting operators of conserved quantities. It is seen that the theta-type invariant spaces for the $BC_N$ Ruijsenaars model correspond to the spaces which are related to the quasi-exact solvability for the $BC_N$ Inozemtsev model, and the degeneration of the $BC_N$ Inozemtsev model (especially for its quasi-exact solvability) is clarified.

In papers \cite{FGGRZ1,FGGRZ2,Tan}, several models which are related to the Inozemtsev model are studied. It would be interesting to link their results with ours.
In \cite{KT,Tak2}, the method of perturbation for the elliptic Calogero-Moser-Sutherland models from the trigonometric models is introduced.
For the Hamiltonian of the $BC_N$ Inozemtsev model, holomorphy of perturbation from the trigonometric model can be established.
Relationship between the perturbation and the complete integrability should be clarified. More precisely, holomorphy of perturbation for commuting operators of conserved quantities should be shown, although it is not successful as of now.

\appendix
\section{} \label{sect:append}
This appendix presents the definitions of and formulae for elliptic functions.

Let $\omega_1$ and $\omega_3$ be complex numbers such that the value $\omega_3/ \omega_1$ is an element of the upper half plane.
The Weierstrass $\wp$-function is defined as follows:
\begin{align}
& \wp (x)=\wp(x|2\omega_1, 2\omega_3)= \label{def:wp} \\
& \; \; \; \; \frac{1}{x^2}+\sum_{(m,n)\in \Zint \times \Zint \setminus \{ (0,0)\} } \left( \frac{1}{(x-2m\omega_1 -2n\omega_3)^2}-\frac{1}{(2m\omega_1 +2n\omega_3)^2}\right),.\nonumber 
\end{align}
Setting $\omega_2=-\omega_1-\omega_3$ and 
\begin{align}
& e_i=\wp(\omega_i) \; \; \; \; \; \; (i=1,2,3) .
\end{align}
yields the relations
\begin{align}
&  e_1+e_2+e_3=0, \; \; \; \wp(x+2\omega_j)=\wp(x); (j=1,2,3), \label{wpformula} \\
& \frac{\wp''(x)}{\wp'(x)^2}=\frac{1}{2}\left( \frac{1}{\wp(x)-e_1}+\frac{1}{\wp(x)-e_2}+\frac{1}{\wp(x)-e_3} \right), \nonumber \\
& \wp (x+y)= \frac{1}{4}\left( \frac{\wp'(x)+\wp'(y)}{\wp(x)-\wp(y)}\right) ^2 -\wp(x)-\wp(y), \nonumber \\
& \wp(x+y) +\wp(x-y) =\frac{\wp'(x)^2+\wp'(y)^2}{2(\wp(x)-\wp(y))^2}-2\wp(x)-2\wp(y), \nonumber \\
& \wp(x+\omega_i)=e_i+\frac{(e_i-e_{i'})(e_i-e_{i''})}{\wp(x)-e_i}, \; \; \; \; (i=1,2,3), \nonumber
\end{align}
where $i', i'' \in \{1,2,3\}$ with $i'<i''$, $i\neq i'$, and $i\neq i''$.

Let $\omega_1=1/2$ and $\tau = \omega_3/ \omega _1$.
The Jacobi theta functions are defined by
\begin{align}
& \theta _1(x)= 2 \sum_{n=1}^{\infty} (-1)^{n-1}e^{\tau \pi \sqrt{-1} (n-\frac{1}{2})^2} \sin (2n-1) \pi x , \label{def:theta} \\
& \theta _2(x)= 2 \sum_{n=1}^{\infty} e^{\tau \pi \sqrt{-1} (n-\frac{1}{2})^2} \cos (2n-1) \pi x  ,\nonumber \\
& \theta _3(x)= 1+2 \sum_{n=1}^{\infty} e^{\tau \pi \sqrt{-1} n^2} \cos 2n \pi x  ,\nonumber \\
& \theta _0(x)= \theta _4(x)= 1+2 \sum_{n=1}^{\infty} (-1)^n e^{\tau \pi \sqrt{-1} n^2} \cos 2n \pi x .\nonumber 
\end{align}
Then the following relations are satisfied:
\begin{align}
& \wp(x)=- \frac{d^2}{dx^2}\log \theta _1(x) + \mbox{const}, \label{thetaperiod}\\
& \wp (x+\omega_i) = - \frac{d^2}{dx^2}\log \theta _{i+1}(x) +\mbox{const} \; \; \; \; (i=1,2,3), \nonumber \\
& \theta _1 (x)=-\theta _1(-x) , \; \; \; \; \theta _i (x)=\theta _i(-x) \; (i=0,2,3), \nonumber \\
& \theta _i (x+1)=-\theta _i(x) \; (i=1,2), \; \; \; \; \theta _i (x+1)=\theta _i(x) \; (i=0,3), \nonumber \\
& \theta _i (x+\tau)=-e^{-\pi \sqrt{-1} (2x+\tau)}\theta _i(x)  \; (i=0,1) , \nonumber \\
& \theta _i (x+\tau)=e^{-\pi \sqrt{-1} (2x+\tau)}\theta _i(x) \; (i=2,3), \nonumber \\
& \theta _1 (2x) \theta _2(0) \theta _3(0) \theta _0(0)= 2\theta _1 (x) \theta _2(x) \theta _3(x) \theta _0(x), \nonumber \\
& \theta ' _1 (0)=\pi \theta _2(0) \theta _3(0) \theta _0(0). \nonumber 
\end{align}

Let $p=\exp (\pi \sqrt{-1} \tau )$. The expansions of the functions $\wp(x)$, $\wp(x+\frac{1}{2} )$, $\wp (x+\frac{\tau}{2} )$ and $\wp(x+\frac{1+\tau}{2} )$ in $p$ are given as follows:
\begin{align}
& \wp (x)=
\frac{\pi^2 }{\sin^2 (\pi x)}- \frac{\pi ^2}{3} -8\pi^2 \sum_{n=1}^{\infty} \frac{np ^{2n}}{1-p ^{2n}} (\cos 2n \pi x -1). \label{wpth} \\
& \wp \left(x+\frac{1}{2} \right)=
\frac{\pi^2 }{\cos ^2 (\pi x)}- \frac{\pi ^2}{3} -8\pi^2 \sum_{n=1}^{\infty} \frac{np ^{2n}}{1-p ^{2n}} ((-1)^n \cos 2n \pi x -1), \nonumber \\
&\wp \left(x+\frac{\tau}{2} \right)=- \frac{\pi ^2}{3}
-8\pi^2 \sum_{n=1}^{\infty} np^n\frac{\cos 2\pi n x-p^n}{1-p ^{2n}}, \nonumber \\
&\wp \left(x+\frac{1+\tau}{2} \right)=- \frac{\pi ^2}{3}
-8\pi^2 \sum_{n=1}^{\infty} np^n\frac{(-1)^n\cos 2\pi n x-p^n}{1-p ^{2n}}. \nonumber
\end{align}


\begin{thebibliography}{9999}
\bibitem{vD}
van Diejen, J. F., Integrability of difference Calogero-Moser systems. {\it J. Math. Phys.} {\bf 35} (1994), 2983--3004. 
\bibitem{FGGRZ1}
Finkel F., Gomez-Ullate D., Gonzalez-Lopez A., Rodriguez M. A., and Zhdanov R., $A\sb N$-type Dunkl operators and new spin Calogero-Sutherland models. {\it Commun. Math. Phys.} {\bf 221} (2001), 477--497.
\bibitem{FGGRZ2}
Finkel F., Gomez-Ullate D., Gonzalez-Lopez A., Rodriguez M. A., and Zhdanov R., New spin Calogero-Sutherland models related to $B\sb N$-type Dunkl operators. {\it Nuclear Phys. B} {\bf 613} (2001), 472--496. 
\bibitem{GGR}
Gomez-Ullate D., Gonzalez-Lopez A., and Rodriguez M. A., Exact solutions of a new elliptic Calogero-Sutherland model. {\it Phys. Lett.} {\bf B511} (2001), 112--118. 
\bibitem{HK}
Hikami K. and Komori Y., Diagonalization of the elliptic Ruijsenaars model of type-$BC$. {\it J. Phys. Soc. Japan} {\bf 67} (1998), 4037--4044. 
\bibitem{Ino}
Inozemtsev V. I., Lax representation with spectral parameter on a torus for integrable particle systems. {\it Lett. Math. Phys.} {\bf 17} (1989), 11-17. 
\bibitem{Kom1}
Komori Y., Theta functions associated with affine root systems and the elliptic Ruijsenaars operators. In {\it Progr. Math.} vol.{\bf 191}, pp. 141--162, Birkhauser Boston, MA, 2000. 
\bibitem{Kom2}
Komori Y., Ruijsenaars' commuting difference operators and invariant subspace spanned by theta functions. {\it J. Math. Phys.} {\bf 42} (2001), 4503--4522.
\bibitem{KH}
Komori Y. and Hikami K., Quantum integrability of the generalized elliptic Ruijsenaars models. {\it J. Phys. A} {\bf 30} (1997), 4341--4364.
\bibitem{KT}
Komori Y. and Takemura K., The perturbation of the quantum Calogero-Moser-Sutherland system and related results. {\it Commun. Math. Phys.} {\bf 227} (2002), 93--118.
\bibitem{OOS}
Ochiai H., Oshima T., and Sekiguchi H., Commuting families of symmetric differential operators. {\it Proc. Japan. Acad.} {\bf 70} (1994) 62--66.
\bibitem{OP}
Olshanetsky M. A. and Perelomov A. M., Quantum integrable systems related to Lie algebras. {\it Phys. Rep.} {\bf 94} (1983), 313--404.
\bibitem{O}
Oshima T., Completely integrable systems with a symmetry in coordinates. {\it Asian J. Math.} {\bf 2} (1998), 935--955. 
\bibitem{OS}
Oshima T. and Sekiguchi H., Commuting families of differential operators invariant under the action of a Weyl group. {\it J. Math. Sci. Univ. Tokyo} {\bf 2} (1995) 1--75.
\bibitem{Rui}
Ruijsenaars S. N. M.: Complete integrability of relativistic Calogero-Moser systems and elliptic function identities. {\it Commun. Math. Phys.} {\bf 110} (1987), 191--213. 
\bibitem{ST}
Sasaki R. and Takasaki K., Quantum Inozemtsev model, quasi-exact solvability and $\mathcal N$-fold supersymmetry. {\it J. Phys. A} {\bf 34} (2001), 9533--9553.
\bibitem{Tak1}
Takemura K., The Heun equation and the Calogero-Moser-Sutherland system I: the Bethe Ansatz method. Preprint, math.CA/0103077, 2001.
\bibitem{Tak2}
Takemura K., The Heun equation and the Calogero-Moser-Sutherland system II: the perturbation and the algebraic solution. Preprint math.CA/0112179, 2001.
\bibitem{Tak3}
Takemura K., The Heun equation and the Calogero-Moser-Sutherland system III: the finite gap property and the monodromy. Preprint math.CA/0201208, 2002.
\bibitem{Tan}
Tanaka T., A Family of Quasi-solvable Quantum Many-body Systems. Preprint hep-th/0202101, 2002.
\bibitem{Tur}
Turbiner A. V., Quasi-exactly-solvable problems and ${\rm sl}(2)$ algebra. {\it Commun. Math. Phys.} {\bf 118} (1988), 467--474.
\end{thebibliography}
\end{document}